\documentclass[12pt,a4paper,twoside,reqno]{article}
\usepackage[utf8]{inputenc}
\usepackage{amsmath}
\usepackage{amsthm}
\usepackage{amssymb}
\usepackage{amsfonts}
\usepackage{microtype}
\usepackage{breakurl}
\usepackage{calrsfs}
\usepackage[usenames]{xcolor}
\usepackage[colorlinks=true,
linkcolor=webgreen,
filecolor=webbrown,
citecolor=webgreen]{hyperref}

\definecolor{webgreen}{rgb}{0,.5,0}
\definecolor{webbrown}{rgb}{.6,0,0}

\usepackage{fullpage}
\usepackage{float}
\numberwithin{equation}{section}
\numberwithin{figure}{section}

\theoremstyle{plain}
\newtheorem{theorem}{Theorem}
\newtheorem{corollary}[theorem]{Corollary}
\newtheorem{lemma}[theorem]{Lemma}
\newtheorem{proposition}[theorem]{Proposition}
\theoremstyle{definition}

\newtheorem{example}[theorem]{Example}

\theoremstyle{remark}
\newtheorem{remark}[theorem]{Remark}

\newcommand{\seqnum}[1]{\href{https://oeis.org/#1}{\rm \underline{#1}}}
\let\cal=\mathcal
\global\long\def\tU{\varLambda}
\date{}

\begin{document}

\title{Sequences Derived from\\
The Symmetric Powers of $\{1,2,\ldots,k\}$}
\author{Po-Yi Huang and Wen-Fong Ke\medskip{}
\\
{\normalsize Department of Mathematics}\\
{\normalsize National Cheng Kung University}\\
{\normalsize Tainan 701, Taiwan}\\
{\normalsize\href{mailto:pyhuang@ncku.edu.tw}{pyhuang@ncku.edu.tw}}\\
{\normalsize\href{mailto:wfke@mail.ncku.edu.tw}{wfke@mail.ncku.edu.tw}}}
\maketitle
\begin{abstract}
For a fixed integer $k$, we define a sequence $A_k=(a_k(n))_{n\geq0}$ and a corresponding sparse subsequence $S_k$ using the cardinality of the $n$-th symmetric power of the set $\{1,2,\ldots, k\}$. For $k\in\{2,\dots,8\}$, we find recursive formulas for $S_k$, and show that the values $a_{k}(0)$, $a_{k}(1)$, and $a_{k}(3)$ are sufficient for constructing $A_{k}$.
\end{abstract}

\section{Introduction}

Let $C$ and $D$ be finite subsets of $\mathbb{N}$. We recall that
the \textit{symmetric difference} of $C$ and $D$, denoted by $C\mathbin{\triangle}D$,
is the union of $C\setminus D$ and $D\setminus C$. If $C=\{c_{1},\ldots,c_{s}\}$
and $D=\{d_{1},\ldots,d_{t}\}$, both nonempty, we set $C\mathbin{\nabla}D=c_{1}D\mathbin{\triangle}c_{2}D\mathbin{\triangle}\cdots\mathbin{\triangle}c_{s}D$,
where $c_{i}D=\{c_{i}d_{1},\ldots,c_{i}d_{t}\}$ for all~$i$. For
example,
\[
\{1,2,3\}\mathbin{\nabla}\{2,4\}=\{2,4\}\mathbin{\triangle}\{4,8\}\mathbin{\triangle}\{6,12\}=\{2,6,8,12\}.
\]
When either $C$ or $D$ is empty, we put $C\mathbin{\nabla D=\varnothing}$.
This binary operation $\nabla$ on the finite subsets of $\mathbb{N}$,
called the \textit{symmetric product,} is associative, commutative
and distributive over~$\triangle$. We let $F$ denote the set of all
finite subsets of $\mathbb{N}$. Consider the map that takes the $i$-th
prime $p_{i}$ to the $i$-th variable $x_{i}$ of $X=\{x_{1},x_{2},\ldots\}$.
We can extend it in a natural way to a map that takes the natural
number $t=2^{s_{1}}3^{s_{2}}\cdots$ ($s_{i}\geq0$ with finitely many
nonzero~$s_{i}$) to the monomial $m_{t}=x_{1}^{s_{1}}x_{2}^{s_{2}}\cdots$,
then further to a map $f$ that takes the nonempty finite subset $T=\{t_{1},\ldots,t_{\ell}\}$
of $\mathbb{N}$ to the multivariate polynomial $f(T)=m_{t_{1}}+\cdots+m_{t_{\ell}}$.
Set $f(\varnothing)=0$. Then $f:F\to\mathbb{Z}_{2}[X]$ is a ring
isomorphism from $(F,\triangle,\nabla)$ to the unique factorization
domain $(\mathbb{Z}_{2}[X],+,\cdot)$ of multivariate polynomials
on $X=\{x_{i}\mid i\in\mathbb{N}\}$ over the Galois field of order
two.

For $k\in\mathbb{N}$, set $H_{k}=\{1,\ldots,k\}$ and define $H_{k}^{\nabla0}=\{1\}$,
$H_{k}^{\nabla1}=H_{k}$, $H_{k}^{\nabla2}=H_{k}\mathbin{\nabla}H_{k}$,
and $H_{k}^{\nabla n}=H_{k}^{\nabla(n-1)}\mathbin{\nabla}H_{k}$ for
$n>2$. We call $H_{k}^{\nabla n}$ the \textit{$n$-th} \textit{symmetric
power} of $H_{k}$. It is clear that $H_{k}^{\nabla2}=\left\{ 1^{2},\ldots,k^{2}\right\} $.
In general, we have $H_{k}^{\nabla2^{t}}=\left\{ 1^{2^{t}},\ldots,k^{2^{t}}\right\} $
for $t\in\mathbb{N}$. One easily sees that $H_{k}^{\nabla(2n)}=(H_{k}^{\nabla n})^{\nabla2}=\left\{ x^{2}\mid x\in H_{k}^{\nabla n}\right\} $
for all~$n\in\mathbb{N}$.

In this paper, we consider the sequences $A_{k}=(a_{k}(n))_{n\geq0}$,
where $a_{k}(n)=|H_{k}^{\nabla n}|$, the cardinality of the
$n$-th symmetric power of $H_{k}$. Thus, for all $k,t\in\mathbb{N}$, we have
$a_{k}(2^{t})=|H_{k}^{\nabla2^{t}}|=k$. These sequences arise
during the investigation of the minimal distances of certain linear
codes (cf.~\cite{CaicedoCP,HuangKP,PachS,Pilz92}) about which we will not
go into detail.

The first few of these sequences can be found in the database \textit{The
On-Line Encyclopedia of Integer Sequences (OEIS) }\cite{OEIS}:
\begin{itemize}
\item $A_{1}=[1,1,1,\ldots]$ is \seqnum{A000012};
\item $A_{2}=[1,2,2,4,2,4,4,8,2,4,4,8,4,8,8,16,2,\ldots]$ is \seqnum{A001316};
\item $A_{3}=[1,3,3,9,3,9,9,27,3,9,9,27,9,27,27,81,3,\ldots]$ is \seqnum{A048883};
\item $A_{4}=[1,4,4,12,4,16,12,40,4,16,16,48,12,48,40,128,4,\ldots]$ is
\seqnum{A253064}.
\end{itemize}
These appear as the numbers of ON cells at the $n$-th generation
when certain ``odd-rule'' cellular automata, with rules 001, 003,
013 and 017, on the square grid
are started in generation $0$ with a single ON cell at the origin (see \cite{Sloane2015}).

The sequences $A_{k}$, $k\geq5$, are not present in the OEIS database.
While it is easy to imagine that $a_{k}(n)$ could get as large as
possible with $n$ when $k\geq2$, the situation that $a_{k}(2^{t})=k$
for all $t\geq0$ makes the sequences $A_{k}$ interesting.

We discuss in this article how one can compute $A_{k}$ when $2\leq k\leq8$.
We show that $a_{k}(n)$ is closely related to the binary expansion
of $n$. Here, if $n=\epsilon_{0}+2\epsilon_{1}+\cdots+2^{\ell}\epsilon_{\ell}$,
where $\ell\geq0$, and each $\epsilon_{i}\in\{0,1\}$ with $\epsilon_{\ell}=1$,
then the binary expansion of $n$ is $[\epsilon_{\ell}\epsilon_{\ell-1}\cdots\epsilon_{1}\epsilon_{0}]$,
and we write $(n)_{2}=[\epsilon_{\ell}\epsilon_{\ell-1}\cdots\epsilon_{1}\epsilon_{0}]$.

In the next section, we make an observation that allows us to restrict
our attention to the \textit{sparse subsequences} $S_{k}=(a_{k}(2^{n}-1))_{n\geq0}$
when $k\leq7$. In Section~\ref{sec:Sparse}, we focus on determining
the recursive formula for $S_{k}$ when $2\leq k\leq8$. In the last
section, we return to the sequence~$A_{8}$.

\section{An observation}\label{sec:Observation}

Note that every number in $H_{k}^{\nabla n}$, $2\leq k\leq8$ and
$0\leq n$, is a product of some powers of $2$, $3$, $5$, and $7$.
\begin{lemma}
\label{lem:1}Let $k,n\in\mathbb{N}$ with $1\leq k\leq7$. If $n=\alpha+\beta\cdot2^{s+1}$
for some $\alpha,\beta,s\in\mathbb{N}$ with $\alpha<2^{s}$, then
$|H_{k}^{\nabla n}|=|H_{k}^{\nabla\alpha}|\cdot|H_{k}^{\nabla\beta}|$.
\end{lemma}

\begin{proof}
We have $H_{k}^{\nabla n}=H_{k}^{\nabla\alpha}\mathbin{\nabla}H_{k}^{\nabla(\beta\cdot2^{s+1})}$.
Let $u\in H_{k}^{\nabla\alpha}$ and $v\in H_{k}^{\nabla(\beta\cdot2^{s+1})}$.
Then there are $e_i,f_i\in \mathbb N\cup\{0\}$, $1\leq i\leq 4$, such that
$ u=2^{e_{1}}3^{e_{2}}5^{e_{3}}7^{e_{4}}$ and $v=2^{f_{1}}3^{f_{2}}5^{f_{3}}7^{f_{4}}$.
Since $k\leq7$ and $\alpha<2^{s}$, we have $e_{i}\leq2\alpha<2^{s+1}$
for each $i$. On the other hand, we have $2^{s+1}\mid f_{i}$ for each $i$.
Hence, we have $H_{k}^{\nabla\alpha}\cap H_{k}^{\nabla(\beta\cdot2^{s+1})}=\varnothing$.
Suppose further that $u'=2^{e_{1}'}3^{e_{2}'}5^{e_{3}'}7^{e_{4}'}\in H_{k}^{\nabla\alpha}$
and $v'=2^{f_{1}'}3^{f_{2}'}5^{f_{3}'}7^{f_{4}'}\in H_{k}^{\nabla(\beta\cdot2^{s+1})}$
are such that $uv=u'v'$. Then, for each $i$, we have $e_{i}+f_{i}=e_{i}'+f_{i}'$, which is equivalent to $f_{i}'-f_{i}=e_{i}-e_{i}'$. Since $2^{s+1}\mid(f_{i}'-f_{i})$
and $|e_{i}-e_{i}'|<2^{s+1}$, we get $f_{i}'-f_{i}=0$. It follows that
$e_{i}=e_{i}'$, and so $u=u'$ and $v=v'$. Therefore, we have $|H_{k}^{\nabla n}|=|H_{k}^{\nabla\alpha}|\cdot|H_{k}^{\nabla\beta}|$
as claimed.
\end{proof}
In terms of binary expansions, the numbers $\alpha$, $\beta$, and
$n$ in Lemma~\ref{lem:1} can be presented as $(\alpha)_{2}=[x\cdots y]$,
$(\beta)_{2}=[b\cdots d]$, and $(n)_{2}=[(b\cdots d)(0\cdots0)(x\cdots y)]$,
where at least one~$0$ is present in $(0\cdots0)$. Thus, in the
cases that $k\leq7$, the lemma implies that if we can compute the
sparse subsequence $S_{k}=(a_{k}(2^{n}-1))_{n\geq0}$, then
we can compute every term $a_{k}(n)$ of $A_{k}$ according to the
binary expansion of $n$. For example, if $(n)_{2}=[1011011101111]$,
that is $n=11727$, then $|H_{k}^{\nabla n}|=|H_{k}^{\nabla1}|\cdot|H_{k}^{\nabla3}|\cdot|H_{k}^{\nabla7}|\cdot|H_{k}^{\nabla15}|$
for $k\in\{1,2,\ldots,7\}$.

In general, Lemma~\ref{lem:1} does not hold for $k\geq8$. For example, we have
$\left|H_{8}^{\nabla3}\right|=48$, $|H_{8}^{\nabla8}|=|H_{8}|=8$,
and $|H_{8}^{\nabla11}|=368\neq8\cdot48$. We will deal with
$A_{8}$ separately in the last section.

\section{The sparse subsequences}\label{sec:Sparse}

In this section, we consider the sparse subsequences $S_{k}=(a_{k}(2^{n}-1))_{n\geq0}$
where $2\leq k\leq 8$. For simplicity, with $k$ fixed, we
write $\theta_{n}=a_{k}(2^{n}-1)=|H_{k}^{\nabla(2^{n}-1)}|$
for $n\geq0$. We have the following recursive formulas for $S_{k}$.

\begin{proposition}
\label{prop:sparse}For each $k$, we have $\theta_{0}=1$ and $\theta_{1}=k$.
\begin{enumerate}
\item If $k=2$, then $\theta_{n+1}=2\theta_{n}$ for $n\geq0$.
\item If $k=3$, then $\theta_{n+1}=3\theta_{n}$ for $n\geq0$.
\item If $k=4$, then $\theta_{n+2}=2\theta_{n+1}+4\theta_{n}$ for $n\geq0$.
\item If $k=5$, then $\theta_{n+2}=3\theta_{n+1}+6\theta_{n}$ for $n\geq0$.
\item If $k=6$, then $\theta_{n+1}=5\theta_{n}$ for $n\geq1$.
\item If $k=7$, then $\theta_{n+2}=6\theta_{n+1}+\theta_{n}$ for $n\geq0$.
\item If $k=8$, then $\theta_{n+3}=7\theta_{n+2}-2\theta_{n+1}-24\theta_{n}$
for $n\geq0$ with $\theta_{2}=48$.
\end{enumerate}
\end{proposition}

Together with Lemma~\ref{lem:1} and the reduction rules to be derived
in Section~\ref{sec:k=00003D8} for $k=8$, this proposition gives
us the following corollary.
\begin{corollary}
\label{cor:reductions}Let $n,k\in\mathbb{N}$.
\begin{enumerate}
\item If $2\leq k\leq7$, then $a_{k}(n)$ can be computed from $a_{k}(0)=1$
and $a_{k}(1)=k$.
\item If $k=8$, then $a_{k}(n)$ can be computed from $a_{k}(0)=1$, $a_{k}(1)=8$,
and $a_{k}(3)=48$.
\end{enumerate}
\end{corollary}

The rest of the section is devoted to showing the validity of the
formulas in Proposition~\ref{prop:sparse}. We start with $S_{8}$
as it provides the best illustration of our approach.

\subsection{The case \texorpdfstring{$k=8$}{k=8}}

For $n\geq0$, the set $H_{8}^{\nabla(2^{n}-1)}$ is the disjoint union of
sets of ``chains''. These chains are classified as types A, B, and~C. A chain of type~A is of the form $\{x,2x,4x,\ldots\}$; a chain of
type~B is of the form $\{x,4x,16x,\ldots\}$; and a chain of type~C
is just a singleton that is not contained in any chain of type~A
or type~B. These chains are uniquely determined via the following
procedure. First, collect all possible chains of type~A. After removing
the elements of these collected chains from $H_{8}^{\nabla(2^{n}-1)}$,
collect all possible chains of type~B from the remaining elements.
After removing the elements of the type~B chains thus collected,
each of the remaining elements is itself a chain of type~C. Let ${\cal C}_{n}$
be the collection of the chains of $H_{8}^{\nabla(2^{n}-1)}$.

For $n=0$, we are looking at $H_{8}^{\nabla(2^{0}-1)}=H_{8}^{\nabla0}=\{1\}$.
The chains in $\mathcal{C}_{0}$ are
\begin{itemize}
\item type A: none;
\item type B: none;
\item type C: $\{1\}$.
\end{itemize}

For $n=1$, we are looking at $H_{8}^{\nabla(2^{1}-1)}=H_{8}$. The
chains in $\mathcal{C}_{1}$ are
\begin{itemize}
\item type A: $\{1,2,4,8\}\text{ and }\{3,6\}$;
\item type B: none;
\item type C: $\{5\}\text{ and }\{7\}$.
\end{itemize}

For $n=2$, we are looking at
\[
H_{8}^{\nabla(2^{2}-1)}=H_{8}^{\nabla3}=H_{8}^{\nabla2}\mathbin{\nabla}H_{8}=
\{1^{2},2^{2},3^{2},4^{2},5^{2},6^{2},7^{2},8^{2}\}\mathbin{\nabla}\{1,2,3,4,5,6,7,8\}.
\]
The chains in $\mathcal{C}_{2}$ are
\begin{itemize}
\item type A:
\begin{alignat*}{2}
 & \{1,2\}, &  & \{3,6,12,24,48,96,192,384\},\\
 & \{9,18\}, &  & \{25,50,100,200\},\\
 & \{27,54,108,216\}, & \qquad & \{49,98,196,392\},\\
 & \{75,150\}, &  & \{144,288\},\\
 & \{147,294\}, &  & \{256,512\};
\end{alignat*}
\item type B:
\begin{alignat*}{2}
 & \{5,20,80,320\}, & \qquad & \{7,28,112,448\},\\
 & \{45,180\}, &  & \{63,252\};
\end{alignat*}
\item type C:
\[
\{125\},\{175\},\{245\},\{343\}.
\]
\end{itemize}

Before we continue, let us define some numbers. For a fixed $n\geq0$,
we look at ${\cal C}_{n}$, and set the numbers $b_{n}$, $c_{n}$, $u_{n}$, $v_{n}$, and $r_{n}$ as follows:
\begin{itemize}
\item $b_{n}$ is the total number of the elements in the chains of type~A;
\item $c_{n}$ is the number of chains of type~A;
\item $u_{n}$ is the total number of the elements in the chains of type~B;
\item $v_{n}$ is the number of chains of type~B;
\item $r_{n}$ is the total number of elements in the chains of type~C.
\end{itemize}
Thus,
\begin{align*}
 & b_{0}=c_{0}=u_{0}=v_{0}=0\text{ and }r_{0}=1;\\
 & b_{1}=6,c_{1}=2,u_{1}=v_{1}=0\text{ and }r_{1}=2;\\
 & b_{2}=32,c_{2}=10,u_{2}=12,v_{2}=4\text{ and }r_{2}=4.
\end{align*}
We refer to $V_{n}=(b_{n},c_{n},u_{n},v_{n},r_{n})^{\text{t}}$, a
column vector, as the \textit{structural vector} of $H_{8}^{\nabla(2^{n}-1)}$.
They can be obtained inductively as we will do next. Note that we
have $$\theta_{n}=|H_{8}^{\nabla(2^{n}-1)}|=b_{n}+u_{n}+r_{n}.$$

Assume that $n\geq2$ and that we have obtained the collection of
chains ${\cal C}_{n}$ and the structural vector $(b_{n},c_{n},u_{n},v_{n},r_{n})^{\text{t}}$
of $H_{8}^{\nabla(2^{n}-1)}$. We want to get ${\cal C}_{n+1}$ and
the structural vector $(b_{n+1},c_{n+1},u_{n+1},v_{n+1},r_{n+1})^{\text{t}}$
of $H_{8}^{\nabla(2^{n+1}-1)}$. Realizing that
\[
H_{8}^{\nabla(2^{n+1}-1)}=(H_{8}^{\nabla(2^{n}-1)})^{\nabla2}\mathbin{\nabla}H_{8}=\triangle_{C\in{\cal C}_{n},D\in{\cal C}_{1}}(C^{\nabla2}\mathbin{\nabla}D),
\]
we compute all $C^{\nabla2}\mathbin{\nabla}D$ for $C\in{\cal C}_{n}$
and $D\in{\cal C}_{1}$.
\begin{enumerate}
\item Suppose that $x_{i}\cdot\{1,2,4,\ldots,2^{\ell_{i}-1}\}$, $i=1,2,\ldots,c_{n}$,
are the chains of type~A in ${\cal C}_{n}$. Take $x=x_{i}$ and
$\ell=\ell_{i}$ ($1\leq i\leq c_{n}$). Then
\[
\{x,2x,4x,\ldots,2^{\ell-1}x\}^{\nabla2}=x^{2}\cdot\{1,4,16,64,\ldots,2^{2(\ell-1)}\}.
\]

\begin{enumerate}
\item We have
\begin{alignat*}{1}
x^{2}\cdot\text{} & \{1,4,16,64,\ldots,2^{2(\ell-1)}\}\mathbin{\nabla}\{1,2,4,8\}\\
 & =x^{2}\cdot\{1,2^{2\ell},2,2^{2\ell+1}\}=\{x^{2},2x^{2}\}\cup\{2^{2\ell}x^{2},2\cdot2^{2\ell}x^{2}\}.
\end{alignat*}
Since $\ell\geq2$, this yields two distinct chains of type~A with
length $2$. Thus, these symmetric products contribute $4c_{n}$ to
$b_{n+1}$, and $2c_{n}$ to $c_{n+1}$.
\item We have
\begin{alignat*}{1}
x^{2}\cdot\text{} & \{1,2,4,\ldots,2^{\ell-1}\}^{\nabla2}\mathbin{\nabla}\{3,6\}\\
 & =3x^{2}\cdot\{1,4,16,\ldots,2^{2(\ell-1)}\}\mathbin{\nabla}\{1,2\}\\
 & =3x^{2}\cdot\{1,2,4,8,\ldots,2^{2\ell-1}\},
\end{alignat*}
a chain of type~A with length $2\ell$. Thus, these contribute $2b_{n}$
to $b_{n+1}$, and $c_{n}$ to~$c_{n+1}$.
\item We have
\[
x^{2}\cdot\{1,2,4,\ldots,2^{\ell-1}\}^{\nabla2}\mathbin{\nabla}\{5\}=5x^{2}\cdot\{1,4,16,\ldots,2^{2(\ell-1)}\}
\]
 and
\[
x^{2}\cdot\{1,2,4,\ldots,2^{\ell-1}\}^{\nabla2}\mathbin{\nabla}\{7\}=7x^{2}\cdot\{1,4,16,\ldots,2^{2(\ell-1)}\}.
\]
Each of them is a chain of type~B with length $\ell$. Such symmetric
products contribute $2b_{n}$ to $u_{n+1}$, and $2c_{n}$ to~$v_{n+1}$.
\end{enumerate}
\item Suppose that $y_{i}\cdot\{1,4,16,\ldots,4^{\ell_{i}'-1}\}$, $i=1,2,\ldots,v_{n}$,
are the chains of type~B in ${\cal C}_{n}$. Take $y=y_{i}$ and
$\ell'=\ell_{i}'$ ($1\leq i\leq v_{n}$). Then
\[
\{y,4y,16y,\ldots,4^{\ell'-1}y\}^{\nabla2}=y^{2}\cdot\{1,16,256,\ldots,4^{2(\ell'-1)}\}.
\]

\begin{enumerate}
\item As
\begin{align*}
\{y,4y,16y,\ldots,4^{\ell'-1}y\}^{\nabla2}\mathbin{\nabla}\{1,2,4,8\} & =cup_{i=0}^{\ell'-1}4^{2i}y^{2}\cdot\{1,2,4,8\}\\
 & =y^{2}\cdot\{1,2,4,8,16,32,\ldots,2^{4\ell'-1}\},
\end{align*}
which is a chains of type~A of length $4\ell'$, these contribute
$4u_{n}$ to $b_{n+1}$, and $v_{n}$ to~$c_{n+1}$.
\item As
\[
\{y,4y,16y,\ldots,4^{\ell'-1}y\}^{2}\mathbin{\nabla}\{3,6\}=cup_{i=0}^{\ell'-1}3\cdot4^{2i}y^{2}\cdot\{1,2\},
\]
which yields $\ell'$ chains of type A each of length $2$, these
contribute $2u_{n}$ to $b_{n+1}$, and $u_{n}$ to~$c_{n+1}$.
\item Each
\[
\{y,4y,16y,\ldots,4^{\ell'-1}y\}^{2}\mathbin{\nabla}\{5\}=cup_{i=0}^{\ell'-1}5\cdot4^{2i}y^{2}\cdot\{1\}
\]
 and
\[
\{y,4y,16y,\ldots,4^{\ell'-1}y\}^{2}\mathbin{\nabla}\{7\}=cup_{i=0}^{\ell'-1}7\cdot4^{2i}y^{2}\cdot\{1\}
\]
yields $2\ell'$ chains of type~C. Thus, these contribute $2u_{n}$
to~$r_{n+1}$.
\end{enumerate}
\item Suppose that $\{z_{i}\}$, $i=1,2,\ldots,r_{n}$, are the chains of
type~C in ${\cal C}_{n}$. Take $z=z_{i}$ ($1\leq i\leq r_{n}$).
\begin{enumerate}
\item Each $\{z\}^{\nabla2}\mathbin{\nabla}\{1,2,4,8\}=\{z^{2},2z^{2},4z^{2},8z^{2}\}$
is a chain of type~A with length $4$. Thus, these contribute $4r_{n}$
to $b_{n+1}$, and $r_{n}$ to $c_{n+1}$.
\item Each $\{z\}^{\nabla2}\mathbin{\nabla}\{3,6\}=\{3z^{2},6z^{2}\}$ gives
a chain of type~A with length $2$. Thus, these contribute~$2r_{n}$
to $b_{n+1}$, and $r_{n}$ to~$c_{n+1}$.
\item Each of $\{z\}^{\nabla2}\mathbin{\nabla}\{5\}=\{5z^{2}\}$ and $\{z\}^{\nabla2}\mathbin{\nabla}\{7\}=\{7z^{2}\}$
gives a chain of type~C. Thus, these contribute~$2r_{n}$ to~$r_{n+1}$.
\end{enumerate}
\end{enumerate}

Summarizing, we get the following system of recursive
linear equations expressing $b_{n+1}$, $c_{n+1}$, $u_{n+1}$, $v_{n+1}$, and $r_{n+1}$
in terms of $b_{n}$, $c_{n}$, $u_{n}$, $v_{n}$, and $r_{n}$:
\[
\left\{\begin{array}{llrrrlr}
b_{n+1}&= & 2b_{n} & +4c_{n} & +6u_{n} &  & +6r_{n},\\
c_{n+1}&= &  & 3c_{n} & +u_{n}\phantom{6} & +v_{n} & +2r_{n},\\
u_{n+1}&= & 2b_{n}& & & &,\\
v_{n+1}&= &  & 2c_{n}& & &,\\
r_{n+1}&= &  &  & 2u_{n} &  & +2r_{n}.
\end{array}\right.\label{eq:linear_recursive}
\]

Now, we need to clarify that the symmetric products of chains as described
above do indeed give us the correct structural vector $V_{n+1}$ of
$H_{8}^{\nabla(2^{n+1}-1)}$ from that of $H_{8}^{\nabla(2^{n}-1)}$.
Specifically, we need to demonstrate the following. Let $C_{1}$ and
$C_{2}$ be two chains in $H_{8}^{\nabla(2^{n}-1)}$, and let $D_{1}$
and $D_{2}$ be two chains in $H_{8}$. If either $C_{1}\not=C_{2}$
or $D_{1}\not=D_{2}$, then
\begin{enumerate}
\item the chains obtained in the symmetric products $C=C_{1}^{\nabla2}\mathbin{\nabla}D_{1}$
and $C'=C_{2}^{\nabla2}\mathbin{\nabla}D_{2}$ are distinct chains
in $H_{8}^{\nabla(2^{n+1}-1)}$, and
\item no longer chains can be formed from these chains.
\end{enumerate}

Since $H_{8}^{\nabla(2^{0}-1)}=H_{8}^{\nabla0}=\{1\}$,
it is clear that the chains in $H_{8}^{\nabla(2^{1}-1)}=H_{8}$ are
exactly produced in this way. Also, one can easily check that the
chains in $H_{8}^{\nabla(2^{2}-1)}=H^{\nabla3}$ can be produced exactly
in this way from those in $H_{8}$.

Now, if we can show that $C\cap C'$ is empty, then (1) is true. Further,
if we can show that $C\cap2C'$, $C'\cap2C$, $C\cap4C'$, and $C'\cap4C$
are all empty, then no concatenating chains to form longer chains
are possible, and so (2) holds.

We will show that these are indeed true after we make the following
general observation.
\begin{lemma}
\label{lem:ind}Let $C_{1}$ and $C_{2}$ be two nonempty finite subsets
of $\mathbb{N}$, and $\ell\geq0$. Then the following statements are
equivalent.
\begin{enumerate}
\item For $j\in\{0,\ldots,\ell\}$, the intersections $C_{1}\cap(2^{j}C_{2})$ and $C_{2}\cap(2^{j}C_{1})$ are empty.
\item If $2^{g_{1}}v_{1}\in C_{1}$ and $2^{g_{2}}v_{2}\in C_{2}$, where
$g_{i}\geq0$ and $2\nmid v_{i}$ for $i=1$ and $2$, then either
$v_{1}\not=v_{2}$ or $|g_{1}-g_{2}|\geq\ell+1$.
\end{enumerate}
\end{lemma}

\begin{proof}
Let $2^{g_{1}}v_{1}\in C_{1}$ and $2^{g_{2}}v_{2}\in C_{2}$, where
$g_{i}\geq0$ and $2\nmid v_{i}$ for $i=1$ and $2$.

Assume (1) and $v_{1}=v_{2}$. If $j\geq0$, then $g_{1}=g_{2}+j$
implies $2^{g_{1}}v_{1}=2^{j}\cdot(2^{g_{2}}v_{2})\in C_{1}\cap(2^{j}C_{2})$,
and $g_{2}=g_{1}+j$ implies $2^{g_{2}}v_{2}=2^{j}\cdot(2^{g_{1}}v_{1})\in C_{2}\cap(2^{j}C_{1})$.
By the assumption, these cannot happen for $j\leq\ell$. Hence $|g_{1}-g_{2}|\geq\ell+1$,
and (2) is true.

Conversely, assume (2). If $2^{g_{1}}v_{1}=2^{g_{2}+j}v_{2}\in C_{1}\cap(2^{j}C_{2})$
or $2^{g_{1}+j}v_{1}=2^{g_{2}}v_{2}\in C_{2}\cap(2^{j}C_{1})$,
where $j\geq0$, then $v_{1}=v_{2}$. In such cases, we
have $j=|g_{1}-g_{2}|\geq\ell+1$. Therefore, we get $C_{1}\cap(2^{j}C_{2})=C_{2}\cap(2^{j}C_{1})=\varnothing$
whenever $j\in\{0,1,\ldots,\ell\}$. This is (1).
\end{proof}

We now state and prove the desired statement.
\begin{proposition}
\label{prop:C_C'8}Let $C_{1},C_{2}\in{\cal C}_{n}$, and $D_{1},D_{2}\in{\cal C}_{1}$.
Put $C=C_{1}^{\nabla2}\mathbin{\nabla}D_{1}$ and $C'=C_{2}^{\nabla2}\mathbin{\nabla}D_{2}$.
If either $C_{1}\not=C_{2}$ or $D_{1}\not=D_{2}$, then $C\cap(2^{i}C')=C'\cap(2^{i}C)=\varnothing$
for $i=0,1,2$.
\end{proposition}

\begin{proof}
Suppose first that $D_{1}\not=D_{2}$. Then $(D_{1}\setminus D_{2})\cup(D_{2}\setminus D_{1})$
is not empty, and an odd prime can be found there. Without loss of
generality, we may assume that $3\in(D_{1}\setminus D_{2}$). If $2^{e_{1}}3^{e_{2}}5^{e_{3}}7^{e_{4}}\in C=C_{1}^{\nabla2}\mathbin{\nabla D_{1}}$,
then $e_{2}$ is odd. On the other hand, for $i\in\{0,1,2\}$, if
$2^{e_{1}'}3^{e_{2}'}5^{e_{3}'}7^{e_{4}'}\in2^{i}C'=2^{i}(C_{2}^{\nabla2}\nabla D_{2})$,
then $e_{2}'$ is even. Hence $C\cap(2^{i}C')=\varnothing$.
Similarly, we have $C'\cap(2^{i}C)=\varnothing$ for $i\in\{0,1,2\}$.

Suppose that $D_{1}=D_{2}=D$ and $C_{1}\not=C_{2}$. Consider $u_{1}^{2}d_{1}\in C_{1}^{\nabla2}\mathbin{\nabla}D$
and $u_{2}^{2}d_{2}\in C_{2}^{\nabla2}\mathbin{\nabla}D$, where $d_{1},d_{2}\in D$,
$u_{1}=2^{g_{1}}v_{1}\in C_{1}$, and $u_{2}=2^{g_{2}}v_{2}\in C_{2}$
with $g_{1}\geq0$, $g_{2}\geq0$, and $v_{1}$ and $v_{2}$ being
odd. As distinct chains of $H_{8}^{\nabla(2^{n}-1)},$ $C_{1}$ and
$C_{2}$ satisfy
\begin{equation}
C_{1}\cap(2^{i}C_{2})=\varnothing\text{ and }C_{2}\cap(2^{i}C_{1})=\varnothing\text{ for }i\in\{0,1,2\}.\label{eq:ind}
\end{equation}
Write $d_{1}=2^{w}v$, where $w\geq0$ and $2\nmid v$. Since $D$
is a chain in $H_{8}$ we have $d_{2}=2^{\pm\lambda}d_{1}=2^{\pm\lambda+w}v$
for some $\lambda$ with $0\leq\lambda\leq3$. Therefore,
\[
u_{1}^{2}d_{1}=2^{2g_{1}+w}v_{1}^{2}v\text{ and }u_{2}^{2}d_{2}=2^{2g_{2}\pm\lambda+w}v_{2}^{2}v.
\]
Now, if $v_{1}^{2}v=v_{2}^{2}v$, then $v_{1}=v_{2}$, and so $|g_{1}-g_{2}|\geq3$
by Lemma~\ref{lem:ind}. Hence
\[
|(2g_{1}+w)-(2g_{2}\pm\lambda+w)|\geq2\cdot|g_{1}-g_{2}|-|\lambda|\geq3.
\]
Again, by Lemma~\ref{lem:ind}, we have $C\cap(2^{i}C')=C'\cap(2^{i}C)=\varnothing$
for $i=0$, $1$, and $2$.
\end{proof}
Expressing the linear recursive system (\ref{eq:linear_recursive})
in matrix form, we have $V_{n+1}=MV_{n}$, where
\begin{equation}
M=\begin{pmatrix}2 & 4 & 6 & 0 & 6\\
0 & 3 & 1 & 1 & 2\\
2 & 0 & 0 & 0 & 0\\
0 & 2 & 0 & 0 & 0\\
0 & 0 & 2 & 0 & 2
\end{pmatrix}\text{ and }V_{n}=\begin{pmatrix}b_{n}\\
c_{n}\\
u_{n}\\
v_{n}\\
r_{n}
\end{pmatrix},\text{ }n=0,1,2,\ldots,\text{ with }V_{0}=\begin{pmatrix}0\\
0\\
0\\
0\\
1
\end{pmatrix}.\label{eq:M8}
\end{equation}
Setting $U=\left(1,0,1,0,1\right)$, we have
\[
\theta_{n+1}=b_{n+1}+u_{n+1}+r_{n+1}=UV_{n+1}=UMV_{n}.
\]

Now, the minimal polynomial of $M$ is $X^{2}(X^{3}-7X^{2}+2X+24)$,
and so we have
\[
U(M^{3}-7M^{2}+2M+24I)M^{2}=0.
\]
A direct check shows that
\begin{align*}
UM & =\left(4,4,8,0,8\right),\\
UM^{2} & =\left(24,28,44,4,48\right),\\
UM^{3} & =\left(136,188,268,28,296\right),
\end{align*}
and that
\begin{equation}
UM^{3}-7UM^{2}+2UM+24U=0.\label{eq:M}
\end{equation}
Therefore, for $n\geq0$, we have
\begin{equation}
\theta_{n+3}=UM^{3}V_{n}=(7UM^{2}-2UM-24U)V_{n}=7\theta_{n+2}-2\theta_{n+1}-24\theta_{n}.\label{eq:r8}
\end{equation}

\begin{remark}
\label{rem:tabularize}We can tabulate the above computations:{\small
\[
\begin{array}{cc||c|c|c||ccccc}
 &  & \{1,2,4,8\} & \{3,6\} & \{5\},\{7\} & b_{n+1} & c_{n+1} & u_{n+1} & v_{n+1} & r_{n+1}\\
\hline\hline (b_{n},c_{n}) & \text{A}_{\ell} & 2\text{A}_{2} & \text{A}_{2\ell} & 2\text{B}_{\ell} & 4c_{n}+2b_{n} & 3c_{n} & 2b_{n} & 2c_{n}\\
\hline (u_{n},v_{n}) & \text{B}_{\ell} & \text{A}_{4\ell} & \ell\text{A}_{2} & 2\ell\text{C} & 4u_{n}+2u_{n} & v_{n}+u_{n} &  &  & 2u_{n}\\
\hline (r_{n}) & \text{C} & \text{A}_{4} & \text{A}_{2} & 2\text{C} & 6r_{n} & 2r_{n} &  &  & 2r_{n}
\\\hline\hline \end{array}.
\]
}{\footnotesize}\\
Here, in the table, we use $\text{A}_{\ell}$ and $\text{B}_{\ell}$ to indicate a chain of type~A
and type~B of length $\ell$, respectively; and we use $\text{C}$ to indicate a chain of a singleton. Also, we use $2\text{A}_{\ell}$
to indicate two disjoint chains of type~A with length $\ell$, etc. The
$(b_{n},c_{n})$, $(u_{n},v_{n})$, and $(r_{n})$ in the first column
are to remind the variables used for counting in the respective rows,
and the break-downs are shown in the columns to the right.

For $D\in\{\{1,2,4,8\},\{3,6\},\{5\},\{7\}\}$, the first row records
the pattern for the symmetric products $\text{A}_{\ell}^{\nabla2}\mathbin{\nabla}D$,
the second row records the pattern for the symmetric products
$\text{B}_{\ell}^{\nabla2}\mathbin{\nabla}D$, and the last row records
the pattern for the symmetric products $\text{C}^{\nabla2}\mathbin{\nabla}D$.
The columns of $b_{n+1},c_{n+1},u_{n+1},v_{n+1},r_{n+1}$ show us
the recursive equations we have seen above. For example, the column
of $b_{n+1}$ reads $b_{n+1}=4c_{n}+2b_{n}+4u_{n}+2u_{n}+6r_{n}=2b_{n}+4c_{n}+6u_{n}+6r_{n}$.
The corresponding matrix $M$ of the system of linear equations can
then be written down easily.
\end{remark}

\subsection{The case \texorpdfstring{$k=2,3,4,5,6,7$}{k=2,3,4,5,6,7}}

Recall that $H_{k}^{\nabla(2^{n+1}-1)}=(H_{k}^{\nabla(2^{n}-1)})^{\nabla2}\mathbin{\nabla}H_{k}=\triangle_{i=1}^{k}i(H_{k}^{\nabla(2^{n}-1)})^{\nabla2}$.
Since the elements in $(H_{2}^{\nabla(2^{n}-1)})^{\nabla2}$
and in $(H_{3}^{\nabla(2^{n}-1)})^{\nabla2}$ are squares
for all $n\geq1$, and those in $2(H_{2}^{\nabla(2^{n}-1)})^{\nabla2}$,
$2(H_{3}^{\nabla(2^{n}-1)})^{\nabla2}$, and $3(H_{3}^{\nabla(2^{n}-1)})^{\nabla2}$
are not, we have
\begin{gather*}
(H_{2}^{\nabla(2^{n}-1)})^{\nabla2}\cap(2(H_{2}^{\nabla(2^{n}-1)})^{\nabla2})=\varnothing,\\
(H_{3}^{\nabla(2^{n}-1)})^{\nabla2}\cap(2(H_{3}^{\nabla(2^{n}-1)})^{\nabla2})=\varnothing,\\
(H_{3}^{\nabla(2^{n}-1)})^{\nabla2}\cap(3(H_{3}^{\nabla(2^{n}-1)})^{\nabla2})=\varnothing.
\end{gather*}
Also it is clear that
\[
(2(H_{3}^{\nabla(2^{n}-1)})^{\nabla2})\cap(3(H_{3}^{\nabla(2^{n}-1)})^{\nabla2})=\varnothing.
\]
Therefore, for $k=2$ and $3$, and $n\geq1$, we have
\begin{align}
\theta_{n+1}=|H_{k}^{\nabla(2^{n+1}-1)}| & =|\triangle_{i=1}^{k}i(H_{k}^{\nabla(2^{n}-1)})^{\nabla2}|\nonumber \\
 & =\sum_{i=1}^{k}|i(H_{k}^{\nabla(2^{n}-1)})^{\nabla2}|\nonumber \\
 & =k\cdot|(H_{k}^{\nabla(2^{n}-1)})^{\nabla2}|=|H_{k}^{\nabla(2^{n}-1)}|\cdot|H_{k}|=k\theta_{n}.\label{eq:n=00003D2,3}
\end{align}
As $\theta_{0}=1$, (\ref{eq:n=00003D2,3}) holds for $n\geq0$.

For each $k\in\{4,5,6,7\}$ we analyze $H_{k}$ and compute the symmetric
products $C^{\nabla2}\mathbin{\nabla}D$ of the chains $C$ of $H_{k}^{\nabla(2^{n}-1)}$
and the chains $D$ of $H_{k}$ to produce a table just like we did
for $H_{8}^{\nabla(2^{n+1}-1)}$ in Remark~\ref{rem:tabularize}.
It turns out that in these cases we only need to collect chains of
types~A and~C. Namely, for each $H_{k}^{\nabla(2^{n}-1)}$, $k\in\{4,5,6,7\}$, and $n\in\mathbb N$, we would collect all possible chains of type~A. After removing the elements of these collected chains, each of the remaining elements gives us a chain of type~C. The structural vector of $H_{k}^{\nabla(2^{n}-1)}$ is then set to be
$V_{n}=(b_{n},c_{n},r_{n})^{\text{t}}$, where
$b_{n}$ is the total number of the elements in the chains of type~A, $c_{n}$ is the number of chains of type~A, and
$r_{n}$ is the total number of elements in the chains of type~C.

We also
need to justify that when the chains $C_{1}$ and $C_{2}$ of $H_{k}^{\nabla(2^{n}-1)}$
are distinct and/or the chains $D_{1}$ and $D_{2}$ of $H_{k}$ are
distinct, the chains produced in the symmetric product $C_{1}^{\nabla2}\mathbin{\nabla}D_{1}$
and those produced in $C_{2}^{\nabla2}\mathbin{\nabla}D_{2}$ are
disjoint and cannot be concatenated to form longer chains. The following
proposition similar to Proposition~\ref{prop:C_C'8} provides us
the justification.
\begin{proposition}
\label{prop:C_C'4-7}Let $k\in\{4,5,6,7\}$. Let $C_{1}$ and $C_{2}$
be chains of $H_{k}^{\nabla(2^{n}-1)}$, and $D_{1}$ and $D_{2}$
that of $H_{k}$. Put $C=C_{1}^{\nabla2}\mathbin{\nabla}D_{1}$ and
$C'=C_{2}^{\nabla2}\mathbin{\nabla}D_{2}$. If either $C_{1}\not=C_{2}$
or $D_{1}\not=D_{2}$, then for $i=0$ and $1$, $C\cap(2^{i}C')=C'\cap(2^{i}C)=\varnothing$.
\end{proposition}

To prove Proposition~\ref{prop:C_C'4-7}, one can use the proof of
Proposition~\ref{prop:C_C'8} as a template. The only changes needed
are to adjust the range of $\lambda$ from $0\leq\lambda\leq3$ to
$0\leq\lambda\leq2$ (as $8$ is not present in $H_{k}$), change
the inequality $|g_{1}-g_{2}|\geq3$ to $|g_{1}-g_{2}|\geq2$,
and apply Lemma~\ref{lem:ind} to obtain the desired conclusion from
the inequality
\[
|(2g_{1}+w)-(2g_{2}\pm\lambda+w)|\geq2\cdot|g_{1}-g_{2}|-|\lambda|\geq2.
\]

Thus, for each $k=4,5,6,7$, $M$ is a $3\times3$ matrix. The size
$\theta_{n+1}=|H_{k}^{\nabla(2^{n+1}-1)}|$ satisfies $\theta_{n+1}=UMV_{n}$,
where $U=\left(1,0,1\right)$ and $V_{n}=\left(b_{n},c_{n},r_{n}\right)^{\text{t}}$.
Therefore, to get the recursive formula for $\theta_{k}$, it suffices
to write down the table of the symmetric products of the chains, the
matrix $M$, the minimal polynomial of $M$ and the polynomial $P$
in $M$ with $UP(M)=0$.

\subsubsection{\texorpdfstring{$k=4$}{k=4}}

The table of the symmetric products of the chains reads

{\small
\[
\begin{array}{cc||c|c||ccc}
 &  & \{1,2,4\} & \{3\} & b_{n+1} & c_{n+1} & r_{n+1}\\
\hline\hline (b_{n},c_{n}) & \text{A}_{\ell} & 2A_{2},(\ell-2)\text{C} & \ell\text{C} & 4c_{n} & 2c_{n} & 2b_{n}-2c_{n}\\
\hline (r_{n}) & \text{C} & \text{A}_{3} & \text{C} & 3r_{n} & r_{n} & r_{n}
\\\hline\hline \end{array}.
\]
}\\
The corresponding matrix $M=\begin{pmatrix}0 & 4 & 3\\
0 & 2 & 1\\
2 & -2 & 1
\end{pmatrix}$ has the minimal polynomial
\[
X^{3}-3X^{2}-2X+4=(X-1)(X^{2}-2X-4).
\]
We have $U(M^{2}-2M-4I)=0$, and so for $n\geq0$,
\[
\theta_{n+2}=2\theta_{n+1}+4\theta_{n}.
\]

\subsubsection{\texorpdfstring{$k=5$}{k=5}}

The table of the symmetric products of the chains reads

{\small
\[
\begin{array}{cc||c|c||ccc}
 &  & \{1,2,4\} & \{3\},\{5\} & b_{n+1} & c_{n+1} & r_{n+1}\\
\hline\hline (b_{n},c_{n}) & \text{A}_{\ell} & 2\text{A}_{2},(\ell-2)\text{C} & 2\ell\text{C} & 4c_{n} & 2c_{n} & 3b_{n}-2c_{n}\\
\hline (r_{n}) & \text{C} & \text{A}_{3} & 2\text{C} & 3r_{n} & r_{n} & 2r_{n}
\\\hline\hline \end{array}.
\]
}\\
The corresponding matrix $M=\begin{pmatrix}0 & 4 & 3\\
0 & 2 & 1\\
3 & -2 & 2
\end{pmatrix}$ has the minimal polynomial
\[
X^{3}-4X^{2}-3X+6=(X-1)(X^{2}-3X-6).
\]
We have $U(M^{2}-3M-6I)=0$, and so for $n\geq0$,
\[
\theta_{n+2}=3\theta_{n+1}+6\theta_{n}.
\]

\subsubsection{\texorpdfstring{$k=6$}{k=6}}

The table of the symmetric products of the chains reads

{\small
\[
\begin{array}{cc||c|c|c||ccc}
 &  & \{3,6\} & \{1,2,4\} & \{5\} & b_{n+1} & c_{n+1} & r_{n+1}\\
\hline\hline (b_{n},c_{n}) & \text{A}_{\ell} & \text{A}_{2\ell} & 2\text{A}_{2},(\ell-2)\text{C} & \ell\text{C} & 2b_{n}+4c_{n} & 3c_{n} & 2b_{n}-2c_{n}\\
\hline (r_{n}) & \text{C} & \text{A}_{2} & \text{A}_{3} & \text{C} & 5r_{n} & 2r_{n} & r_{n}
\\\hline\hline \end{array}.
\]
}The corresponding matrix $M=\begin{pmatrix}2 & 4 & 5\\
0 & 3 & 2\\
2 & -2 & 1
\end{pmatrix}$ has the minimal polynomial
\[
X^{3}-6X^{2}+5X=X(X-1)(X-5).
\]
We have $U(M^{2}-5M)=0$, and so for $n\geq1$,
\[
\theta_{n+1}=5\theta_{n}.
\]

\subsubsection{\texorpdfstring{$k=7$}{k=7}}

The table of the symmetric products of the chains reads

{\small
\[
\begin{array}{cc||c|c|c||ccc}
 &  & \{3,6\} & \{1,2,4\} & \{5\},\{7\} & b_{n+1} & c_{n+1} & r_{n+1}\\
\hline\hline (b_{n},c_{n}) & \text{A}_{\ell} & \text{A}_{2\ell} & 2\text{A}_{2},(\ell-2)\text{C} & 2\ell\text{C} & 2b_{n}+4c_{n} & 3c_{n} & 3b_{n}-2c_{n}\\
\hline (r_{n}) & \text{C} & \text{A}_{2} & \text{A}_{3} & 2\text{C} & 5r_{n} & 2r_{n} & 2r_{n}
\\\hline\hline \end{array}.
\]
}The corresponding matrix $M=\begin{pmatrix}2 & 4 & 5\\
0 & 3 & 2\\
3 & -2 & 2
\end{pmatrix}$ has the minimal polynomial
\[
X^{3}-7X^{2}+5X+1=(X-1)(X^{2}-6X-1).
\]
We have $U(M^{2}-6M-I)=0$, and so for $n\geq0$,
\[
\theta_{n+2}=6\theta_{n+1}+\theta_{n}.
\]

\section{The sequence \texorpdfstring{$A_{8}$}{A\_8}}\label{sec:k=00003D8}

Recall that for the sparse subsequence $S_{8}=(\theta_{n})_{n\geq0}$,
where $\theta_{n}=|H_{8}^{\nabla(2^{n}-1)}|$, we associate a structural vector $V_{n}=(b_{n},c_{n},u_{n},v_{n},r_{n})^{\text{t}}$ with each $H_{8}^{\nabla(2^{n}-1)}$,
so that $\theta_{n}=UMV_{n}$ for all $n\geq0$ with the vector $U=(1,0,1,0,1)$
and the $5\times5$ matrix $M$ in (\ref{eq:M8}). We can extend this
idea to $A_{8}$ in order to get every $a_{8}(n)=|H_{8}^{\nabla n}|$,
$n\geq0$.

For each $n$, we again collect the chains of types~A, B and~C in $H_{8}^{\nabla n}$, and define the structural vector $\tU_n$ of $H_8^{\nabla n}$ to be $(b_n,c_n,u_n,v_n,r_n)^{\text{t}}$. Note that the same variables $b_{n}$, $c_{n}$, $u_{n}$, $v_{n}$, and $r_{n}$ are used here. Yet they are not the same in the indices as we have used them when dealing with the sparse subsequences~$S_{k}$.

Now, the structural vector of $H_{8}^{\nabla0}=\{1\}$ is $\tU_{0}=(0,0,0,0,1)^{\text{t}}=V_{0}$. Next, with $H_{8}^{\nabla1}=\left(H_{8}^{\nabla(2^{0}-1)}\right)^{2}\nabla H_{8}$,
we apply the matrix $M$ in (\ref{eq:M8}) to $\tU_{0}$ to get $\tU_{1}$
(which is $V_{1}$ in this case). Thus, $\tU_{1}=M\tU_{0}=(6,2,0,0,2)^{\text{t}}$.
To get $\tU_{2}$, the structural vector of $H_{8}^{\nabla2}$, we
need another matrix for \textit{symmetric squaring}.

Let $n\geq1$. Suppose that we have obtained the collection of chains $\cal C$ of $H_{8}^{\nabla n}$ and the structural vector $\tU_{n}=(b_n,c_n,u_n,v_n,r_n)^{\text{t}}$. For a type~A chain of length $\ell$, say\break $C=\{x,2x,\ldots,2^\ell x\}$, we have $C^{\nabla 2}=\{x^2,4x^2,\ldots,4^{\ell-1} x^2\}$, which is type~B chain of length~$\ell$ in $H_8^{\nabla 2n}$. For a type~B chain of length $\ell$, say $C=\{x,4x,\ldots,4^{\ell-1} x\}$, we have $C^{\nabla 2}=\{x^2,16x^2,\ldots,16^{\ell-1} x^2\}=\cup_{i=0}^{\ell-1} \{16^ix^2\}$, which is the union of $\ell$ type~C chains in $H_8^{\nabla 2n}$. Obviously, $\{x\}^{\nabla 2}=\{x^2\}$ for every type~C chain $\{x\}$ in $\cal C$. This is a type C chain in $H_8^{\nabla 2n}$. Thus, the table for
symmetric squaring is as follows:{\small
\[
\begin{array}{cc||c||ccccc}
 & C & C^{\nabla2} & b_{2n} & c_{2n} & u_{2n} & v_{2n} & r_{2n}\\
\hline\hline (b_{n},c_{n}) & \text{A}_{\ell} & \text{B}_{\ell} &  &  & b_{n} & c_{n}\\
\hline (u_{n},v_{n}) & \text{B}_{\ell} & \ell\text{C} &  &  &  &  & u_{n}\\
\hline (r_{n}) & \text{C} & \text{C} &  &  &  &  & r_{n}
\\\hline\hline \end{array}.
\]
}The corresponding matrix $W$ is
\[
W=\begin{pmatrix}0 & 0 & 0 & 0 & 0\\
0 & 0 & 0 & 0 & 0\\
1 & 0 & 0 & 0 & 0\\
0 & 1 & 0 & 0 & 0\\
0 & 0 & 1 & 0 & 1
\end{pmatrix}.
\]
We have $\tU_{2n}=W\tU_{n}$ and, in particular, $\tU_{2}=W\tU_{1}=(WM)\tU_{0}$.

\begin{lemma}
\label{lem:WM}Let $n\geq0$ be arbitrary, and suppose that $(n)_{2}=[\epsilon_{\ell}\epsilon_{\ell-1}\cdots\epsilon_{1}\epsilon_{0}]$.
Then
\[
\tU_{n}=E_{0}E_{1}\cdots E_{\ell-1}E_{\ell}\tU_{0},
\]
where, for each $i$, $E_{i}=M$ if $\epsilon_{i}=1$, and $E_{i}=W$
if $\epsilon_{i}=0$.
\end{lemma}

\begin{proof}
The result holds for $\ell=0$ and $\ell=1$ as we have just seen. Let $\ell\geq 2$.
Take $m\in\mathbb{N}$ with $(m)_{2}=[\epsilon_{\ell}\epsilon_{\ell-1}\cdots\epsilon_{1}]$. Assume that $\tU_{m}=E_{1}\cdots E_{\ell-1}E_{\ell}\tU_{0}$ with $E_{i}=M$ if $\epsilon_{i}=1$, and $E_{i}=W$
if $\epsilon_{i}=0$.

If $\epsilon_{0}=0$, then $n=2m$. In this case, we have $H_{8}^{\nabla n}=H_{8}^{\nabla(2m)}=\left(H_{8}^{\nabla m}\right)^{\nabla2}$
and so $\tU_{n}=W\tU_{m}$.

If $\epsilon_{0}=1$, then $n=2m+1$. In this case, we have
\[
H_{8}^{\nabla n}=H_{8}^{\nabla(2m+1)}=H_{8}^{\nabla(2m)}\mathbin{\nabla}H_{8}=\left(H_{8}^{\nabla m}\right)^{\nabla2}\mathbin{\nabla}H_{8},
\]
and so $\tU_{n}=M\tU_{m}$.

Therefore, by the induction hypothesis, we have
\[
\tU_{n}=E_{0}E_{1}\cdots E_{\ell-1}E_{\ell}\tU_{0}
\]
where $E_{i}=M$ if $\epsilon_{i}=1$, and $E_{i}=W$ if $\epsilon_{i}=0$.
\end{proof}
Recall that $U=(1,0,1,0,1)$, and $|H_{8}^{\nabla n}|=U\tU_{n}$ for all $n\geq0$.
We shall derive some reduction rules to help to compute every $|H_{8}^{\nabla n}|$,
$n\in\mathbb{N}$. We already know that
\begin{equation}
|H_{8}^{\nabla2m}|=|H_{8}^{\nabla m}|\text{ for all }m\in\mathbb{N}.\label{eq:8_0}
\end{equation}

Notice that
\[
W^{2}=\tU_{0}U,
\]
and we have the following result similar to Lemma~\ref{lem:1}.
\begin{lemma}
\label{lem:R00}If $n=\alpha+\beta\cdot2^{s+2}$ for some $\alpha,\beta,s\in\mathbb{N}$
with $\alpha<2^{s}$, then
\begin{equation}
|H_{8}^{\nabla n}|=|H_{8}^{\nabla\alpha}|\cdot|H_{8}^{\nabla\beta}|.\label{eq:8_1}
\end{equation}
\end{lemma}

\begin{proof}
Let $(\alpha)_{2}=[x\cdots z]$ and $(\beta)_{2}=[b\cdots d]$. Then
$(n)_{2}=[(b\cdots d)(00\cdots0)(x\cdots z)]$, where there are at least
two $0$'s between $(b\cdots d)$ and $(x\cdots z)$. Let $C_{1}$ and
$C_{2}$ be the products of $M$'s and $W$'s such that $\tU_{\alpha}=C_{1}\tU_{0}$
and $\tU_{\beta}=C_{2}\tU_{0}$. Then for some $\ell\geq0$, we have
\begin{align*}
|H_{8}^{\nabla n}|=U\tU_{n} & =UC_{1}W^{\ell+2}C_{2}\tU_{0}\\
 & =UC_{1}\tU_{0}UW^{\ell}C_{2}\tU_{0}\\
 & =(U\tU_{\alpha})(U\tU_{2^{\ell}\beta})\\
 & =|H_{8}^{\nabla\alpha}|\cdot|H_{8}^{\nabla(2^{\ell}\beta)}|\\
 & =|H_{8}^{\nabla\alpha}|\cdot|(H_{8}^{\nabla\beta})^{\nabla2^{\ell}}|=|H_{8}^{\nabla\alpha}|\cdot|H_{8}^{\nabla\beta}|.
\end{align*}
\end{proof}

Next, we collect some useful facts about $M$ and $W$.
\begin{lemma}
\label{lem:WMeqs}With $M$, $W$, $U$, and $\tU_{0}$ as above, we
have
\begin{enumerate}
\item $U\tU_{0}=1$;
\item $UW=U$;
\item $UMW=8U$;
\item $UM^{2}W=UM+40U$.
\end{enumerate}
Consequently, if $C$ is a product of $M$'s and $W$'s, then $UC=U(pM^{2}+qM+rI)$
for some $p,q,r\in\mathbb{Z}$.
\end{lemma}

\begin{proof}
The four identities can be checked directly. Combining (2)--(4) and
the fact that $UM^{3}=7UM^{2}-2UM-24U$ (see (\ref{eq:M})), we infer that
$UC=U(pM^{2}+qM+rI)$ for some $p,q,r\in\mathbb{Z}$ as claimed.
\end{proof}

\begin{lemma}
\label{lem:R01}If $(n)_{2}=[(b\cdots d)01]$ and $(m)_{2}=[b\cdots d]$,
then
\begin{equation}
|H_{8}^{\nabla n}|=8\cdot|H_{8}^{\nabla m}|.\label{eq:8_2}
\end{equation}
\end{lemma}

\begin{proof}
As $\tU_{n}=MW\tU_{m}$, we have, by Lemma~\ref{lem:WMeqs}(3),
\[
|H_{8}^{\nabla n}|=U\tU_{n}=UMW\tU_{m}=8U\tU_{m}=8\cdot|H_{8}^{\nabla m}|.
\]
\end{proof}
\begin{lemma}
\label{lem:R01-1}If $(n)_{2}=[(b\cdots d)011]$, $(m_{1})_{2}=[b\cdots d1]$
and $(m_{2})_{2}=[b\cdots d]$, then
\begin{equation}
|H_{8}^{\nabla n}|=|H_{8}^{\nabla m_{1}}|+40\cdot|H_{8}^{\nabla m_{2}}|.\label{eq:8_3}
\end{equation}
\end{lemma}

\begin{proof}
As $\tU_{n}=M^{2}W\tU_{m_{2}}$, we have, by Lemma~\ref{lem:WMeqs}(4),
\begin{align*}
|H_{8}^{\nabla n}|=U\tU_{n} & =UM^{2}W\tU_{m_{2}}=(UM+40U)\tU_{m_{2}}\\
 & =UM\tU_{m_{2}}+40U\tU_{m_{2}}=U\tU_{m_{1}}+40U\tU_{m_{2}}=|H_{8}^{\nabla m_{1}}|+40\cdot|H_{8}^{\nabla m_{2}}|.
\end{align*}
\end{proof}

\begin{lemma}
\label{lem:R111}Let $(n)_{2}=[(b\cdots d)111(x\cdots z)]$. With $$(m_{1})_{2}=[(b\cdots d)11(x\cdots z)],\
(m_{2})_{2}=[(b\cdots d)1(x\cdots z)],\text{ and }(m_{3})_{2}=[(b\cdots d)(x\cdots z)],$$
we have
\begin{equation}
|H_{8}^{\nabla n}|=7\cdot|H_{8}^{\nabla m_{1}}|-2\cdot|H_{8}^{\nabla m_{2}}|-24\cdot|H_{8}^{\nabla m_{3}}|,\label{eq:8-4}
\end{equation}
\end{lemma}

\begin{proof}
Let $\alpha\geq0$ and $\beta\geq0$ be such that $(\alpha)_{2}=[x\cdots z]$
and $(\beta)_{2}=[b\cdots d]$, and let $C$ be the products of $M$'s
and $W$'s such that $\tU_{\alpha}=C\tU_{0}$. Then $\tU_{n}=CM^{3}\tU_{\beta}$,
$\tU_{m_{1}}=CM^{2}\tU_{\beta}$, $\tU_{m_{2}}=CM\tU_{\beta}$, and
$\tU_{m_{3}}=C\tU_{\beta}$. By Lemma~\ref{lem:WMeqs}, we have $UC=U(pM^{2}+qM+rI)$
for some $p,q,r\in\mathbb{Z}$. Let $D=pM^{2}+qM+rI$. From (\ref{eq:M}),
we get
\begin{align*}
|H_{8}^{\nabla n}|=U\tU_{n} & =UCM^{3}\tU_{\beta}\\
 & =UDM^{3}\tU_{\beta}\\
 & =UM^{3}D\tU_{\beta}\\
 & =U(7M^{2}-2M-24I)D\tU_{\beta}\\
 & =7UM^{2}D\tU_{\beta}-2UMD\tU_{\beta}-24UD\tU_{\beta}\\
 & =7UDM^{2}\tU_{\beta}-2UDM\tU_{\beta}-24UD\tU_{\beta}\\
 & =7UCM^{2}\tU_{\beta}-2UCM\tU_{\beta}-24UC\tU_{\beta}\\
 & =7U\tU_{m_{1}}-2U\tU_{m_{2}}-24U\tU_{m_{3}}\\
 & =7\cdot|H_{8}^{\nabla m_{1}}|-2\cdot|H_{8}^{\nabla m_{2}}|-24\cdot|H_{8}^{\nabla m_{3}}|.
\end{align*}
\end{proof}

We are in a position to show that we can get every $|H_{8}^{\nabla n}|$,
$n\in\mathbb{N}$, with just $|H_{8}^{\nabla0}|$, $|H_{8}^{\nabla1}|$
and $|H_{8}^{\nabla3}|$ by applying (\ref{eq:8_0})--(\ref{eq:8_3})
repeatedly.

Let ${\cal S}$ be the set of all $n$, $n\geq0$, such that none
of the rules (\ref{eq:8_0}), (\ref{eq:8_1}), (\ref{eq:8_2}) and
(\ref{eq:8_3}) can apply. In particular, if $n\in{\cal S}$, then
the binary expansion $(n)_{2}$, when read from the higher bits to
the lower bits, does not have $00$ or $111$ anywhere, and does not
end with $0$ or~$01$.
\begin{itemize}
\item If $(n)_{2}$ has just one bit, then $(n)_{2}$ is either $(0)_{2}=[0]$, or $(1)_{2}=[1]$.
We have $|H_{8}^{\nabla0}|=1$ and $|H_{8}^{\nabla1}|=8$.
\item If $(n)_{2}$ has two bits, then $(n)_{2}$  has to be $(3)_{2}=[11]$. We have
$|H_{8}^{\nabla3}|=48$.
\end{itemize}
So we have $\{0,1,3\}\subseteq{\cal S}$.
\begin{proposition}
\label{prop:S}The set ${\cal S}$ is exactly $\{0,1,3\}$.
\end{proposition}

\begin{proof}
As $(4)_{2}=[100]$, $(5)_{2}=[101]$, $(6)_{2}=[110]$, and $(7)_{2}=[111]$,
no integers with a $3$ bit binary expansion belong to ${\cal S}$.

Let $n\in{\cal S}$. Suppose that $(n)_{2}$ has at least $4$ bits.
Then the last three bits of $(n)_{2}$ has to be~$011$. So $(n)_{2}=[(b\cdots d)011]$.
But then (\ref{eq:8_3}) applies, contradicting to that $n\in{\cal S}$.
Therefore, $(n)_{2}$ has at most $3$ bits, and we are done.
\end{proof}
Now Corollary~\ref{cor:reductions}(2) follows.

\begin{example}
\begin{enumerate}
  \item The identity (\ref{eq:8_3}) gives us
\[
\left|H_{8}^{\nabla11}\right|=\left|H_{8}^{\nabla3}\right|+40\cdot\left|H_{8}^{\nabla1}\right|=368
\]
since $(11)_{2}=[1011]$, $(3)_{2}=[11]$, and $(1)_{2}=[1]$.
  \item Using (\ref{eq:8_3}) again, we get
\[
\left|H_{8}^{\nabla27}\right|=\left|H_{8}^{\nabla7}\right|+40\cdot\left|H_{8}^{\nabla3}\right|
=\left|H_{8}^{\nabla7}\right|+1920
\]
since $(27)_{2}=[11011]$, $(7)_{2}=[111]$, and $(3)_{2}=[11]$.
It follows from (\ref{eq:8-4}) that
\[
\left|H_{8}^{\nabla7}\right|=7\cdot\left|H_{8}^{\nabla3}\right|-2\cdot\left|H_{8}^{\nabla1}\right|-24\cdot\left|H_{8}^{\nabla0}\right|=296.
\]
Hence $\left|H_{8}^{\nabla27}\right|=2216$.
\end{enumerate}
\end{example}

We give two more reduction rules which can be useful for computing
$|H_{8}^{\nabla n}|$. Here is the first one.
\begin{lemma}
\label{lem:R011}If $(n)_{2}=[(b\cdots d)011011]$, $(m_{1})_{2}=[(b\cdots d)011]$
and $(m_{2})_{2}=[b\cdots d]$, then
\begin{equation}
|H_{8}^{\nabla n}|=47\cdot|H_{8}^{\nabla m_{1}}|-40\cdot|H_{8}^{\nabla m_{2}}|.\label{eq:8_4}
\end{equation}
\end{lemma}

\begin{proof}
We have
\begin{align*}
|H_{8}^{\nabla n}| & =UM^{2}WM^{2}W\tU_{m_{2}}\\
 & =(UM+40U)M^{2}W\tU_{m_{2}}\\
 & =UM^{3}W\tU_{m_{2}}+40UM^{2}W\tU_{m_{2}}\\
 & =(7UM^{2}-2UM-24UI)W\tU_{m_{2}}+40UM^{2}W\tU_{m_{2}}\\
 & =47UM^{2}W\tU_{m_{2}}-2UMW\tU_{m_{2}}-24UW\tU_{m_{2}}\\
 & =47U\tU_{m_{1}}-16U\tU_{m_{2}}-24U\tU_{m_{2}}\\
 & =47U\tU_{m_{1}}-40U\tU_{m_{2}}\\
 & =47\cdot|H_{8}^{\nabla m_{1}}|-40\cdot|H_{8}^{\nabla m_{2}}|.
\end{align*}
\end{proof}
The second one requires the following facts about $WM$, which can
be checked directly. We note that the minimal polynomial of $WM$
is $X^{4}-9X^{3}+8X^{2}$, and so
\begin{equation}
UWM((WM)^{2}-9WM+8I)=0.\label{eq:UWM}
\end{equation}
Since $UW=U$, (\ref{eq:UWM}) is reduced to
\begin{equation}
UM((WM)^{2}-9WM+8I)=0.\label{eq:UWM-1}
\end{equation}
Furthermore, it is also true that
\begin{equation}
UM^{2}((WM)^{2}-9WM+8I)=0.\label{eq:UWM-2}
\end{equation}
Although $UM^{3}((WM)^{2}-9WM+8I)\not=0$, we still have
\begin{equation}
UM^{3}((WM)^{2}-9WM+8I)\tU_{0}=0.\label{eq:UMW-3}
\end{equation}

\begin{lemma}
\label{lem:R1010}Let $n,m_{1},m_{2}\in\mathbb{N}$. If one of the
following holds
\begin{enumerate}
\item $(n)_{2}=[(b\cdots d)101011]$, $(m_{1})_{2}=[(b\cdots d)1011]$ and
$(m_{2})_{2}=[(b\cdots d)11]$,
\item $(n)_{2}=[10101(x\cdots z)]$, $(m_{1})_{2}=[101(x\cdots z)]$, and $(m_{2})_{2}=[1(x\cdots z)]$,
\end{enumerate}
then
\begin{equation}
|H_{8}^{\nabla n}|=9\cdot|H_{8}^{\nabla m_{1}}|-8\cdot|H_{8}^{\nabla m_{2}}|.\label{eq:8_5}
\end{equation}

\end{lemma}

\begin{proof}
(1) Let $\alpha\geq0$ with $(\alpha)_{2}=[b\cdots d]$. We have
\begin{align*}
  |H_{8}^{\nabla m_{2}}| & =UM^{2}\tU_{\alpha}, \\
  |H_{8}^{\nabla m_{1}}| & =UM^{2}WM\tU_{\alpha}, \\\noalign{\hbox{and}}
  |H_{8}^{\nabla n}| & =UM^{2}(WM)^{2}\tU_{\alpha}.
\end{align*}
Together with (\ref{eq:UWM-2}), this gives
\[
|H_{8}^{\nabla n}|-9\cdot|H_{8}^{\nabla m_{1}}|+8\cdot|H_{8}^{\nabla m_{2}}|=UM^{2}((WM)^{2}-9WM+8I))\tU_{\alpha}=0.
\]
(2) Let $\beta\geq0$ with $(\beta)_{2}=[x\cdots z]$. Let $C$ be
the product of $M$'s and $W$'s such that $\tU_{\beta}=C\tU_{0}$.
We have $|H_{8}^{\nabla m_{2}}|=UCM\tU_{0}$, $|H_{8}^{\nabla m_{1}}|=UCMWM\tU_{0}$
and $|H_{8}^{\nabla n}|=UCM(WM)^{2}\tU_{0}$. Let $p,q,r\in\mathbb{Z}$
such that $UC=pUM^{2}+qUM+rU$. Then
\begin{align*}
 & \hspace{-1cm}|H_{8}^{\nabla n}|-9\cdot|H_{8}^{\nabla m_{1}}|+8\cdot|H_{8}^{\nabla m_{2}}|\\
 & =UCM((WM)^{2}-9WM+8I)\tU_{0}\\
 & =U(pM^{2}+qM+rI)M((WM)^{2}-9WM+8I)\tU_{0}\\
 & =pUM^{3}((WM)^{2}-9WM+8I)\tU_{0}\\
 & \qquad+qUM^{2}((WM)^{2}-9WM+8I)\tU_{0}\\
 & \qquad\quad+rUM((WM)^{2}-9WM+8I)\tU_{0}.
\end{align*}
According to (\ref{eq:UWM-1})--(\ref{eq:UMW-3}), the last three
terms are all zero. We are done.
\end{proof}
Here is an example of $|H_{8}^{\nabla n}|$ with a bit larger
$n$.
\begin{example}
\label{exa:9}Let us compute $|H_{8}^{\nabla n}|$ where $(n)_{2}=(11101011011)$,
so $n=1883$.

First of all, we apply (\ref{eq:8-4}) with $(m_{1})_{2}=[11101011]$
and $(m_{2})_{2}=[11101]$ to get
\[
|H_{8}^{\nabla n}|=47\cdot|H_{8}^{\nabla m_{1}}|-40\cdot|H_{8}^{\nabla m_{2}}|,
\]
Next, we put together (\ref{eq:8_5}) and (\ref{eq:8_2}) with $(\beta_{1})_{2}=[111011]$,
$(\beta_{2})_{2}=[1111]$, and $(\gamma)_{2}=[111]$, and we get
\[
|H_{8}^{\nabla m_{1}}|=9\cdot|H_{8}^{\nabla\beta_{1}}|-8\cdot|H_{8}^{\nabla\beta_{2}}|\text{ and }|H_{8}^{\nabla m_{2}}|=8\cdot|H_{8}^{\nabla\gamma}|.
\]
Lastly, we infer from (\ref{eq:8-4}) that
\begin{align*}
|H_{8}^{\nabla\beta_{1}}| & =7\cdot|H_{8}^{\nabla27}|-2\cdot|H_{8}^{\nabla11}|-24\cdot|H_{8}^{\nabla3}|=13624,\\
|H_{8}^{\nabla\beta_{2}}| & =7\cdot|H_{8}^{\nabla7}|-2\cdot|H_{8}^{\nabla3}|-24\cdot|H_{8}^{\nabla1}|=1784,\\
|H_{8}^{\nabla\gamma}| & =|H_{8}^{\nabla7}|=296,
\end{align*}
and so $|H_{8}^{\nabla m_{1}}|=108344$ and $|H_{8}^{\nabla m_{2}}|=2368$.
Therefore, $|H_{8}^{\nabla n}|=4997448$.
\end{example}

When applying the above rules for computing $|H_{8}^{\nabla n}|$,
there is no particular order that one have to follow. Some orders
could be more efficient than others when applied in specific cases.
Still, one can design certain general algorithm to compute $|H_{8}^{\nabla n}|$.

Finally, we remark that we do not have any formula for $A_{k}$ and
$S_{k}$ when $k\geq9$. The approach using chains like above fails
since new types of chains need to be introduced, e.g.,~chains of the
type $\{x,3x,9x,\ldots\}$, and they may not be disjoint from the chains
of other types.

\section{Acknowledgments}

The authors like to thank the reviewer and the editor for their careful reading and comments to improve the clarity and the readability of the paper.

\bigskip
\hrule
\bigskip

\noindent 2020 {\it Mathematics Subject Classification: }Primary 11B37; Secondary 11Y55.

\noindent \emph{Keywords: }symmetric product,
symmetric power, binary expansion, odd-rule cellular automaton.

\bigskip
\hrule
\bigskip

\noindent (Concerned with sequences
\seqnum{A000012},
\seqnum{A001316},
\seqnum{A048883}, and
\seqnum{A253064}.)

\bigskip
\hrule
\bigskip
\end{document}